\documentclass[10pt]{article}   
\usepackage{amssymb} 
\usepackage{latexsym} 
\usepackage{amsthm} 
\usepackage[dvips]{epsfig}

\newcommand{\N}{\mathbb N}                                                      
\newcommand{\Z}{ \mathbb Z}

\newcommand{\Sym}{\mathbb S}

\theoremstyle{plain}
\theoremstyle{remark}

\newtheorem{theo}{\bf Theorem}

\newtheorem{cor}{\bf Corollary}
\newtheorem{prop}{\bf Proposition}
\newtheorem{rem}{\bf Remark}

\newenvironment{preuve}{\begin{trivlist}\item{\bf{Proof.}}}%
{\hfill\rule{2mm}{2mm}\end{trivlist}}

\newcommand{\arbre}{\mbox{\Large\it a}} 
                               
\newcommand{\ao}{{\arbre}_{o}} 
 
\newcommand{\aobt}{{\mbox{\Large\it a}}_{o}^{\: \underline{\diamond}}} 
\newcommand{\aot}{{\mbox{\Large\it a}}_{o}^{\: \diamond}}

\newcommand{\tildeaobt}{\widetilde{{\mbox{\Large\it a}}}_{o}^{\: \underline{\diamond\!}}} 
\newcommand{\tildeaot}{\widetilde{{\mbox{\Large\it a}}}_{o}^{\: \diamond}}     
\newcommand{\tildeaonp}{\widetilde{\arbre}_{o}}

\newcommand{\abtpair}{{\mbox{\Large\it a}}^{\: \underline{\diamond}}} 
\newcommand{\atpair}{{\mbox{\Large\it a}}^{\: \diamond}}

\newcommand{\tildeatpair}{\widetilde{{\arbre}}^{\: \diamond}} 
\newcommand{\tildeabtpair}{\widetilde{{\mbox{\Large\it a}}}^{\: \underline{\diamond}}} 
\newcommand{\tildeaotpair}{\widetilde{{\arbre}}_{o}^{\: \diamond}} 
\newcommand{\tildeaobtpair}{\widetilde{{\mbox{\Large\it a}}}_{o}^{\: \underline{\diamond}}}

\newcommand{\anp}{\arbre} 
\newcommand{\ab}{\arbre^{-}} 
\newcommand{\abt}{{\arbre}^{\: \underline{\diamond}}} 
\newcommand{\at}{{\arbre}^{\: \diamond}} 
\newcommand{\asym}{{\arbre}_{\rm TS}}

\newcommand{\tr}{\:\!\diamond\!} 
\newcommand{\btr}{\:\underline{\!\diamond\!}}

\begin{document} 
\title{Labelled and unlabelled enumeration\\ 
of $k$-gonal 2-trees\footnote{With the support of FCAR (Qu\'ebec) and NSERC (Canada).}}
\author{Gilbert Labelle, C\'edric Lamathe and Pierre Leroux} 
\date{December 22, 2003}
\maketitle
\begin{abstract}
In this paper\footnote{This is the full version of a paper presented
at the International Colloquium on ``Mathematics and Computer Science`` held in
Versailles, France, in September 2002 (see \cite{LLLV}).}, we generalize %
2-trees by replacing triangles by quadrilaterals, pentagons or
$k$-sided polygons ($k$-gons), where $k\geq 3$ is given. This
generalization, to $k$-gonal 2-trees, is natural and is closely
related, in the planar case, to some specializations of the
cell-growth problem. Our goal is the labelled and unlabelled
enumeration of $k$-gonal 2-trees according to the number $n$ of
$k$-gons. We give explicit formulas in the labelled case, and, in
the unlabelled case, recursive and asymptotic formulas.
\end{abstract}
%
\section{Introduction} 
\begin{figure}[h] 
 \centerline{\includegraphics[width=.75\textwidth]{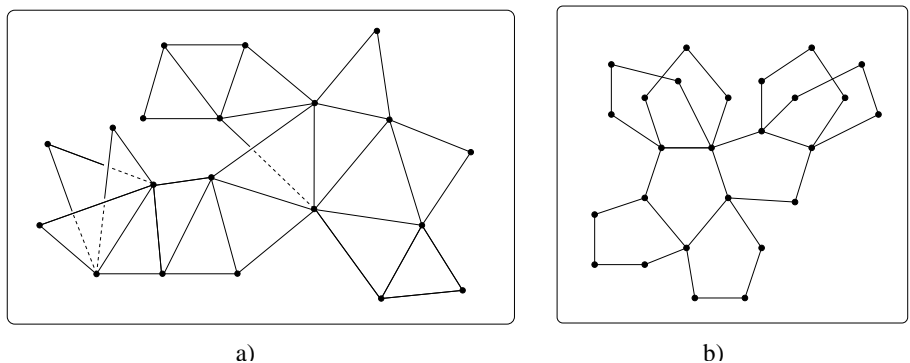}} 
  \caption{$k$-gonal 2-trees with $k=3$ and $k=5$} 
  \label{fig:ajout} 
\end{figure}

\begin{figure}[h] 
 \centerline{\includegraphics[width=.80\textwidth]{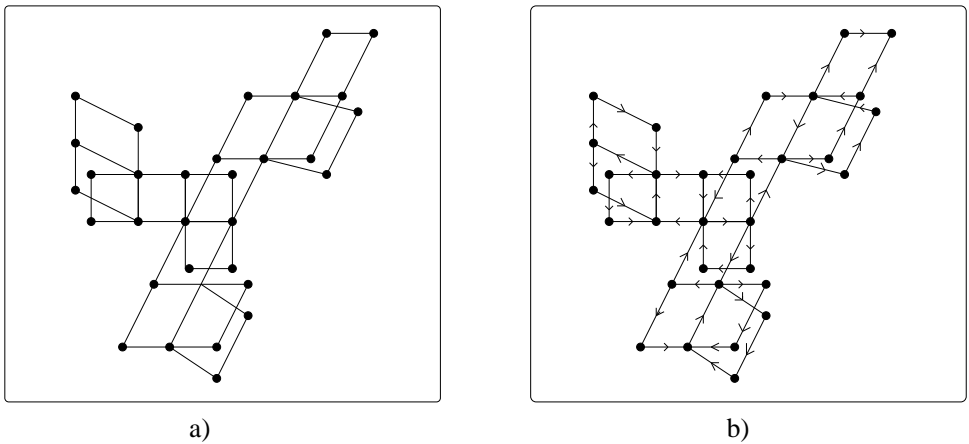}} 
  \caption{Unoriented and oriented 4-gonal 2-trees} 
  \label{fig:exemple} 
\end{figure}

Essentially, a 2-\emph{tree} (or \emph{bidimensional tree}) is a connected
simple graph composed of triangles glued along their edges in a
tree-like fashion, that is, without cycles (of triangles). 
This definition can be extended by replacing the triangles %
by quadrilaterals, pentagons or $k$-sided polygons ($k$-gons), where
$k\geq 3$ is fixed. Such 2-trees, built on $k$-gons, are called {\it$k$-gonal 2-trees}.  %
Figures~\ref{fig:ajout}a, \ref{fig:ajout}b, and \ref{fig:exemple}a
show examples of $k$-gonal 2-trees, for $k=3,5$ and 4, respectively. %
Of course the usual 2-trees correspond to $k=3$.

The enumeration of 2-trees is extensively studied in the literature.  %
The first results in this direction are found in 1970, in Palmer \cite{Pa} %
for the labelled enumeration of 2-trees (see also Beineke and Moon \cite{BM}) %
and in Harary and Palmer \cite{HP} (1973) for the unlabelled enumeration. %
During the same period, Palmer and Read \cite{PaR} enumerated  %
labelled and unlabelled \emph{outerplanar} 2-trees, that is, 2-trees %
which can be embedded in the plane in such a way that each vertex belongs %
to the external face.  %
The term \emph{planar} is also used in this sense. %
See also Labelle, Lamathe and Leroux \cite{LLL1,LLL2}. %
 
Two years later, together with Harary, these authors generalized their results %
in \cite{HPR} by considering for the first time $k$-gonal 2-trees %
and enumerating them in the outerplanar case, in the context of a cell-growth %
problem.  %
 
In his 1993 Ph.D. Thesis \cite{kloksdutch,kloksthese}, %
Ton Kloks enumerated unlabelled {\it biconnected partial 2-trees}, %
that is, 2-trees in which some edges have been deleted without however %
losing the 2-connectedness. He calls these graphs \emph{2-partials}. %
This class strictly contains that of $k$-gonal 2-trees %
since, in a 2-partial, polygons of different sizes can occur %
and some edges can be missing, provided that they are incident to at least three %
polygons. In principle Kloks' method, which extends the traditional  %
\emph{dissimilarity characteristic} of Otter \cite{Otter} to 2-partials, could be used %
to enumerate $k$-gonal 2-trees (with $k$ fixed). However,  %
to our knowledge, this work has not been done.  %

More recently, in 2000,  %
Fowler, Gessel, Labelle and Leroux \cite{TF1,TF2},  %
have proposed some new functional equations %
for the class of (ordinary) 2-trees, 
which yield recurrences and asymptotic formulas for their unlabelled %
enumeration.  %
Their approach, which is based on the theory of combinatorial species of  %
Joyal (see \cite{Joyal, BLL}), is more structural, replacing a potential %
dissimilarity characteristic formula for each individual 2-tree by a Dissymmetry Theorem  %
for the species of 2-trees. Such a theorem can be formulated for most classes of tree-like %
stuctures, for example ordinary (one-dimensional, Cayley) trees or more generaly %
simple graphs, all of whose 2-connected components are in a given class (see \cite{BLL}), %
plane embedded trees (see \cite{LaLe}), various classes of cacti (see \cite{BBLL}, %
etc.  %
 
In the present paper, we extend to $k$-gonal 2-trees the work of Fowler et als,  %
which corresponds to the case $k=3$.  %
In particular, we label the 2-trees at their $k$-gons. %
Our goal is their labelled and unlabelled enumeration,  %
according to the number of $k$-gons. We will give explicit
formulas in the labelled case and recursive and asymptotic formulas in
the unlabelled case, emphasizing the dependency on $k$.  %
Special attention must be given to the cases where $k$ is even. %

We say that a $k$-gonal 2-tree is {\em oriented} if its edges are
oriented in such a way that each $k$-gon forms an oriented cycle; see
Figure~\ref{fig:exemple} b). In fact, for any $k$-gonal 2-tree $s$,
the orientation of any one of its edges can be extended uniquely to
all of $s$ by first orienting all the polygons to which the edge
belongs and then continuing recursively on all adjacent polygons. The
coherence of the extension is ensured by the arborescent (acyclic)
nature of 2-trees.

We denote by $\arbre$ and $\arbre_o$ the species of
$k$-gonal 2-trees and of oriented $k$-gonal 2-trees. For these
species, we use the symbols $-$, $\diamond$ and $\underline{\diamond}$
as upper indices to indicate that the structures are pointed at an
edge, at a $k$-gon, and at a $k$-gon having itself a distinguished
edge, respectively.
 
A first step is the extension to the $k$-gonal case of the Dissymmetry Theorem  %
for 2-trees, which links together these various pointed species.  %
The proof is similar to the case $k=3$ and is
omitted (see \cite{TF1,TF2}). %
\medskip
\begin{theo}{\textsc{Dissymmetry theorem for $k$-gonal 2-trees.}} 
  The species $\arbre_o$ and $\arbre$ of oriented and unoriented
  $k$-gonal 2-trees, respectively, satisfy the following isomorphisms of
  species:
\begin{eqnarray} 
\arbre_o^{-}+\aot&=&\ao+\aobt,\\ 
\ab+\at&=&\anp+\abt. 
\end{eqnarray}
\end{theo} 

There is yet another species to introduce, which plays an essential
role in the process. It is the species $B=\arbre^{\rightarrow}$ of
oriented-edge rooted ($k$-gonal) 2-trees, that is of 2-trees where an
edge is selected and oriented. As mentionned above, the orientation of
the rooted edge can be extended uniquely to an orientation of the
2-tree so that there is a canonical isomorphism $B=\arbre_o^{-}$.  %
However, it is often useful
not to perform this extension and to consider that only the rooted
edge is oriented.

In the next section, we characterize the species
$B=\arbre^{\rightarrow}$ by a combinatorial functional equation and
give some of its consequences. The goal is then to express the various
pointed species occuring in the Dissymmetry Theorem in terms of $B$
and to deduce enumerative results for the species $\arbre_o$ and
$\arbre$.  The oriented case is simpler and carried out first, in
Section 3. The unoriented case is analyzed in Section 4,  %
where $\arbre$ is viewed as a quotient species of $\arbre_o$ %
and two cases are distinguished, according to the parity of the integer
$k$. Finally, asymptotic results are presented in Section 5.

For our purposes, the main tool of species theory is the P\'olya-Robinson-Joyal
Composition Theorem which can be stated as follows (see \cite{BLL}, Th. 1.4.2):  %
let the species $F$ be the (partitionnal) composition of two species,  %
$F=G\circ H=G(H)$.
Then, the exponential generating function
\[
F(x) = \sum_{n\geq 0}f_n{x^n\over n!},
\]
where $f_n=|F[n]|$ is the number of labelled $F$-structures over a set of  %
cardinality $n$, and the \emph{tilde} generating function
\[
\widetilde{F}(x) = \sum_{n\geq 0}\widetilde{f}_n x^n,
\]
where $\widetilde{f}_n=|F[n]/\Sym_n|$ is the number of unlabelled
$F$-structures of order $n$, satisfy the following equations: %
\begin{eqnarray}
F(x) &=& G(H(x)),\\
\widetilde{F}(x) &=& Z_G(\widetilde{H}(x),\widetilde{H}(x^2),\ldots )\label{eq:**}, %
\end{eqnarray} %
where $Z_G(x_1,x_2,\ldots )$ is the cycle index series of $G$. %
Moreover, we have  %
\begin{equation} %
Z_{F}(x_1,x_2,\ldots) = Z_G\circ Z_H = Z_G(Z_H(x_1,x_2,\ldots),Z_H(x_2,x_4,\ldots),\dots). %
\label{eq:ZG(H)} %
\end{equation} %
Here the operation $\circ$ is the \emph{plethystic composition} of symmetric functions %
when the $x_1,x_2,\dots$ are interpreted as \emph{power sum} symmetric functions %
in some other set of variables ${\mathbf s}=(s_1,s_2,s_3,\dots)$:  
$x_i=p_i=p_i(s_1,s_2,\dots):=\sum_{j\geq1}s_j^i$.  %

This interpretation of the cycle index series as symmetric functions can be taken %
as an alternate definition, as follows (see \cite{BLL}, Example 2.3.15 and Rem. 4.3.8).  %
An $F$-structure is said to be \emph{colored} if the elements of its underlying %
set are assigned colors in the set $\{1,2,3,\dots\}$. Such a colored structure has %
a weight $w$ given  by its color distribution monomial %
in the variables ${\mathbf s}=(s_1,s_2,s_3,\dots)$. %
Let us denote by $F(1_{\mathbf s})$ the weighted set of unlabelled colored $F$-structures. %
Its total weight (or \emph{inventory}) $|F(1_{\mathbf s})|_w$ is a symmetric %
function in the variables ${\mathbf s}$ and thus has a unique expression %
in terms of the power sums $x_i=p_i(s_1,s_2,\dots)$ given precisely by $Z_{F}$: %
\begin{equation} %
|F(1_{\mathbf s})|_w = Z_{F}(x_1,x_2,\ldots). %
\label{eq:ZFalt)} %
\end{equation} %

For example, for the species $E_2$, of 2-element sets, and $E$, of sets, we have %
\begin{equation} %
Z_{E_2}(x_1,x_2,\ldots) =\sum_{i<j}s_is_j + \sum_is_i^2  %
= {1\over2}\left((\sum_is_i)^2  + \sum_is_i^2\right) = {1\over2}(x_1^2 + x_2) %
\label{eq:ZE2} %
\end{equation} %
and %
\begin{equation} %
Z_{E}(x_1,x_2,\ldots) = h(s_1,s_2,\dots) = \exp\left(\sum_{i\geq1}\frac{x_i}{i}\right), %
\label{eq:ZE} %
\end{equation} %
where $h=\sum_{n\geq0}h_n$ denotes the complete homogeneous symmetric function. %

%
\section{The species $B$ of oriented-edge rooted 2-trees}
\begin{figure}[h] 
 \centerline{\includegraphics[width=.8\textwidth]{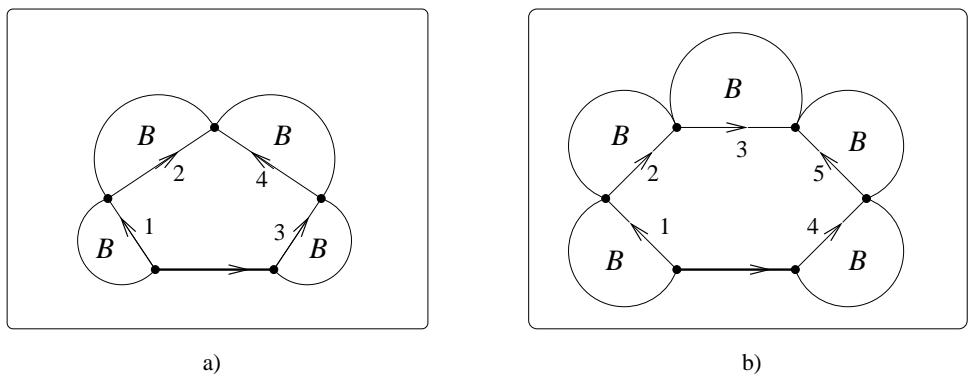}} 
  \caption{A page of an oriented-edge rooted 2-tree, for a) $k=5$, b) $k=6$} 
  \label{fig:page} 
\end{figure}

The species $B=\arbre^{\rightarrow}$ plays a central role in the study
of $k$-gonal 2-trees. The following theorem is an extension to a
general $k$ of the case $k=3$. Note that formula (\ref{eq:aoriente1})
below also makes sense for $k=2$ and corresponds to edge-labelled
(ordinary) rooted trees.
\medskip
\begin{theo}
  The species $B=\arbre^{\rightarrow}$ of oriented-edge rooted
  $k$-gonal 2-trees satisfies the following functional equation
  (isomorphism):
\begin{equation}
B   =  E(XB^{k-1}),\label{eq:aoriente1}
\end{equation}
where $E$ represents the species of sets and $X$ is the species of
singleton $k$-gons.
\end{theo}
\begin{preuve} 
  We decompose an $\arbre^{\rightarrow}$-structure as a set of {\it
    pages}, that is, of maximal subgraphs sharing only one $k$-gon
  with the rooted edge. For each page, the orientation of the rooted
  edge permits to define a linear order and an orientation on the
  $k-1$ remaining edges of the polygon having this edge, in some
  conventional way, for example in the fashion illustrated in
  Figure~\ref{fig:page}a, for the odd case, and \ref{fig:page}b, for the even
  case.  These edges being oriented, we can glue on them some
  $B$-structures.  We then deduce relation (\ref{eq:aoriente1}).
\end{preuve}
%
%

Among the possible edge orientations of an oriented-edge  %
rooted $k$-gon, the one illustrated in Figure \ref{fig:page}a, %
"away from the root edge", has the advantage of remaining valid %
if the root edge is not oriented, for $k$ odd. If $k$ is even, %
we see a difference caused by the existence of an opposite edge %
whose orientation will remain ambiguous. %

We can easily relate the species $B=\arbre^{\rightarrow}$ to that  %
of (ordinary) rooted trees, denoted by $A$, characterized by the
functional equation $A=XE(A)$, where $X$ now represents the sort of
vertices. Indeed from (\ref{eq:aoriente1}), we deduce %
\begin{equation}
(k-1)XB^{k-1} = (k-1)XE((k-1)XB^{k-1}),
\end{equation}
knowing that $E^m(X)=E(mX)$. %
By the Implicit Species Theorem of Joyal (see \cite{BLL}), there %
exists a unique (up to isomorphism) species $Y$ such that $Y=(k-1)XE(Y)$, namely $Y=A((k-1)X)$. %
It follows that %
\begin{equation}
(k-1)XB^{k-1} = A((k-1)X)\label{cest5}
\end{equation} %
and %
\begin{equation}
B^{k-1} = \frac{A((k-1)X)}{(k-1)X}.\label{cest5bis}
\end{equation} %
In analogy with formal power series, it can be shown %
that for any rational number $r\not=0$, any species $F$ %
with constant term equal to 1 (that is $F(0)=1$) admits a \emph{unique} $r^{th}$-root  %
with constant term 1, that is a unique species $G$ such that $G^r=F$ and $G(0)=1$; %
here $G$ may be a virtual species, with rational coefficients (see Rem.~2.6.16 of \cite{BLL}).
In the present case, since both $B$ and $A((k-1)X)/(k-1)X$ have constant term 1, %
we obtain the following expression for the species $B$ in
terms of the species of rooted trees. %
This expression can be used to compute the first terms of the molecular expansion %
of $B$, using Newton's Binomial Theorem; see \cite{ALL03}.
\medskip
\begin{prop} 
  The species $B=\arbre^{\rightarrow}$ of oriented-edge-rooted
  $k$-gonal 2-trees satisfies
\begin{equation} \label{B=F(A)}
B=\sqrt[k-1]{\frac{A((k-1)X)}{(k-1)X}}.\label{eq:B} 
\end{equation} 
\end{prop} 
\medskip
\begin{cor} 
The numbers $a^{\rightarrow}_n$, $a_{n_1,n_2,\ldots }^{\rightarrow}$,
and $b_n=\widetilde{a}_{n}^{\rightarrow}$ of $k$-gonal 2-trees pointed
at an oriented edge and having $n$ $k$-gons, respectively labelled,
fixed by a permutation of cycle type $1^{n_1}2^{n_2}\ldots $ and
unlabelled, satisfy the following formulas and recurrence:
\begin{equation} 
a^{\rightarrow}_n=((k-1)n+1)^{n-1}=m^{n-1},\label{eq:betiqu} 
\end{equation} 
where $m=(k-1)n+1$ is the number of edges,
\begin{equation} 
a_{n_1,n_2,\ldots }^{\rightarrow}=\prod_{i=1}^{\infty}(1+(k-1)\sum_{d|i}dn_d)^{n_i-1}(1+(k-1)\sum_{d|i, d<i}dn_d), 
\label{eq:Zb} 
\end{equation} 
and 
\begin{equation} 
b_n=\frac{1}{n}\sum_{1\leq j\leq n}\sum_{\alpha}(|\alpha|+1)b_{\alpha_1}
b_{\alpha_2}\ldots b_{\alpha_{k-1}}b_{n-j}, \quad \quad b_0=1,\label{eq:bnonetiq} 
\end{equation} 
the last sum running over $(k-1)$-tuples of integers
$\alpha=(\alpha_1,\alpha_2,\ldots ,\alpha_{k-1})$ such that
$|\alpha|+1$ divides the integer $j$, where $|\alpha|=\alpha_1+\alpha_2+\cdots +\alpha_{k-1}$.
\end{cor} 
\begin{preuve} 
Formulas (\ref{eq:betiqu}) and (\ref{eq:Zb}) are obtained by
specializing with $\mu={(k-1)^{-1}}$ the following formulas, given
by Fowler et al. in \cite{TF1,TF2},
\begin{equation}
\left({A(x)\over x}\right)^\mu=\sum_{n\geq 0}\mu(\mu+n)^{n-1}{x^n\over n!},\label{eq:betiq}
\end{equation}
$
Z_{\left({A(X/\mu)\over X/\mu}\right)^\mu}=
$
\begin{equation}
\sum_{n_1,n_2,\ldots}\frac{x_1^{n_1}x_2^{n_2}\ldots}{1^{n_1}n_1!2^{n_2}n_2!\ldots}
\prod_{i=1}^{\infty}(1+{1\over \mu}\sum_{d|i}dn_d)^{n_i-1}(1+{1\over \mu}\sum_{d|i,d<i}dn_d).
\end{equation}
\noindent 
Formula (\ref{eq:betiqu}) can also be established by a Pr\"ufer-like
bijection; see \cite{Prufer, LLLA}. To obtain the recurrence (\ref{eq:bnonetiq}), it suffices
to take the logarithmic derivative of the equation
\begin{equation}\label{eq:btilde}
\widetilde{B}(x)=\exp\left(\sum_{i\geq1}{x^i\widetilde{B}^{k-1}(x^i)\over i}\right)\!, 
\end{equation} 
where $\widetilde{B}(x)=\sum_{n\geq 0}b_nx^n$, which follows from
relation (\ref{eq:aoriente1}), using (\ref{eq:**}) and (\ref{eq:ZE}).
\end{preuve} 

The sequences $\{b_n\}_{n\in \N}$, for $k=2,3,4,5,6$, are listed in
the Encyclopedia of Integer Sequences \cite{Sl,SP}. Respectively: 
A000081, for the number of rooted trees with $n$ nodes,  %
A005750, in relation with planted matched trees with $n$ nodes and 2-trees,  %
A052751, A052773, A052781, in relation with equation (\ref{eq:btilde}).  %
Also, equation (\ref{eq:aoriente1}), is referenced in the Encyclopedia of  %
Combinatorial Structures \cite{INRIA}.  %

Observe that for each $n\geq 1$, $b_n$ is a
polynomial in $k$ of degree $n-1$. This follows from (\ref{eq:Zb}) and
the following explicit formula for $b_n$,
\begin{equation}
b_n = \sum_{n_1+2n_2+\cdots =n}\frac{a^{\rightarrow}_{n_1,n_2,\ldots }}{1^{n_1}n_1!2^{n_2}n_2!\ldots },
\end{equation}
which is a consequence of Burnside's lemma. The asymptotic behavior of
the numbers $b_n$ as $n \rightarrow \infty$, is studied, in particular
as a function of $k$, in Section 7.


%
%
\section{Oriented $k$-gonal 2-trees}  
We begin by determining relations for the pointed species appearing in
the Dissymmetry Theorem. These relations are quite direct and the
proof is left to the reader.
%
%
\begin{prop}\label{prop1}
  The species $\arbre_o^{-}$, $\aot$, and $\aobt$ are characterized by
  the following isomorphisms:
\begin{equation} 
\arbre_o^{-}=B,\quad   
\aot   =  XC_k(B),\quad
\aobt  =  XB^k,    
\end{equation} 
where $B=\arbre^{\rightarrow}$ and $C_k$ represents the species of
oriented cycles of length $k$.
\end{prop} 
 %
Recall that the cycle index series of $C_k$ is given by  %
$Z_{C_k} = \frac{1}{k}\sum_{d|k}\phi(d)x_d^{n/d}$ where $\phi$ is the Euler function.
The Dissymmetry Theorem then permits us to express the ordinary (tilde) generating
series $\widetilde{\arbre}_o(x)$ of unlabelled oriented $k$-gonal
2-trees in terms of the corresponding series for the rooted species:
\begin{equation} 
\tildeaonp(x)  =  \widetilde{\arbre}_o^{-}(x)+\tildeaot(x)-\tildeaobt(x).
\end{equation} 
By Proposition~\ref{prop1}, we can now express
$\widetilde{\arbre}_o(x)$ as function of
$\widetilde{B}(x)=\widetilde{\arbre}^{\rightarrow}(x)$.
%
%
\begin{prop} 
The ordinary generating series $\widetilde{\arbre}_o(x)$ of
unlabelled oriented $k$-gonal 2-trees is given by
\begin{equation}\label{eq:aoriente}
\tildeaonp(x) = \widetilde{B}(x)+{x\over k}\sum_{d|k\atop d>1}\phi(d)\widetilde{B}^{k\over d}(x^d)
-{k-1\over k}x\widetilde{B}^k(x).
\end{equation} 
\end{prop} 
\begin{cor}\label{corol} 
  The numbers $a_{o,n}$ and $\widetilde{a}_{o,n}$ of oriented
  $k$-gonal 2-trees labelled and unlabelled, over $n$ $k$-gons,
  respectively, are given by
\begin{eqnarray} 
a_{o,n} &=& ((k-1)n+1)^{n-2}=m^{n-2},\quad n\geq2,\label{eq:aoetiq}\\ 
\widetilde{a}_{o,n} &=& b_n-\frac{k-1}{k}b_{n-1}^{(k)}+
{1\over k}\sum_{d|k\atop d>1}\phi(d)b_{n-1\over d}^{({k\over d})},\label{eq:aotilden} 
\end{eqnarray}
where 
\[
b_i^{(j)}\ =\ [x^i]\widetilde{B}^j(x)\ =\ \displaystyle{\sum_{i_1+\cdots +i_j=i}b_{i_1}b_{i_2}\ldots b_{i_j}},
\]
denotes the coefficient of $x^i$ in the series
$\widetilde{B}^j(x)$, with $b_r^{(j)}=0$ if $r$ is non-integral or
negative.
\end{cor}  
\begin{preuve} 
For the labelled case, it suffices to remark that
$a_{n}^{\rightarrow}=ma_{o,n}$. In the unlabelled case, equation
(\ref{eq:aotilden}) is directly obtained from (\ref{eq:aoriente}).
\end{preuve} 
 
%
\section{Unoriented $k$-gonal 2-trees}
For the enumeration of (unoriented) $k$-gonal 2-trees, we
consider quotient species of the form $F/\Z_2$, where $F$ is
a species of ``oriented'' structures, %
$\Z_2=\{1,\tau\}$, is a group of order 2  %
and the action of $\tau$ is to reverse the structure orientations.  %
A structure of such a quotient species then consists in an orbit
$\{s,\tau\cdot s\}$ of $F$-structures under the action of $\Z_2$.
 
For instance, the different pointed species of unoriented $k$-gonal
2-trees $\ab$, $\at$ and $\abt$, can be expressed as quotient species
of the corresponding species of oriented $k$-gonal 2-trees:
\begin{equation} 
\ab  = \frac{\arbre^{\rightarrow} }{\Z_2},\quad  
\at  = \frac{\aot}{\Z_2}\ =\ \frac{XC_k(B)}{\Z_2},\quad
\abt = \frac{\aobt}{\Z_2}      \  = \  \frac{XB^k}{\Z_2}. \label{eq:abarre}
\end{equation}
The three basic generating series associated to
such a quotient species, are given by  %
\begin{equation} %
\label{eq:F/Z2(x)} 
(F/\Z_2)(x) = {1\over 2}(F(x) + \sum_{n\geq 0}|{\mathrm{Fix}}_{F_n}(\tau)|{x^n\over n!}), 
\end{equation} %
\begin{equation} %
\label{eq:F/Z2tilde} 
(F/\Z_2)^{\sim}(x)={1\over 2}(\widetilde{F}(x) +  %
  \sum_{n\geq 0}|{\rm{Fix}}_{\widetilde{F}_n}(\tau)|x^n), 
\end{equation} %
where ${\mathrm{Fix}}_{F_n}(\tau)$ and ${\rm{Fix}}_{\widetilde{F}_n}(\tau)$ denote  %
the sets of labelled and unlabelled, respectively, $F$-structures left fixed by the
action of $\tau$, that is, by orientation reversal, and %
\begin{equation} %
\label{eq:ZF/Z2}  %
Z_{F/\Z_2}(x_1,x_2,\dots)= {1\over 2}(Z_F(x_1,x_2,\dots) +  %
 |{\mathrm{Fix}}_{F(1_{\mathbf s})}(\tau)|_w ), %
\end{equation} %
where ${\mathrm{Fix}}_{F(1_{\mathbf s})}(\tau)$ is the set of unlabelled colored %
$F$-structures left fixed by $\tau$, weighted by the color distribution monomials %
in the variables ${\mathbf s}=(s_1,s_2,s_3,\dots)$ and where  %
the inventory $|{\mathrm{Fix}}_{F(1_{\mathbf s})}(\tau)|_w$, being a symmetric function in %
${\mathbf s}$, is expressed in terms of the power sums $x_i=p_i({\mathbf s})$. %
A simple example %
is given by the species $E_2=X^2/\Z_2$, the species of 2-element sets,  %
where formula (\ref{eq:ZF/Z2}) yields immediately $Z_{E_2}={1\over 2}(x_1^2+x_2)$. %

However, some important differences appear in the 
computations, according to the parity of $k$.  %
The main difference comes from the existence of  %
\emph{opposite} edges in $k$-gons, when $k$ is even. %
Accordingly, it is better to treat the two cases separately.
%
%
%
\subsection{Case $k$ odd}  %
If $k$ is odd, 
it is quite simple to extend the method of Fowler et als \cite{TF1,TF2} where $k=3$. %
For example, %
the only labelled oriented $k$-gonal 2-tree left fixed by an orientation reversal,  %
for a given number of polygons, is the one in which all polygons share one common
edge. Hence, from (\ref{eq:F/Z2(x)}) and the fact that $\arbre=\ao/\Z_2$, we deduce %
directly the following.
%
%
\begin{prop}  
If $k$ is odd, the number $a_n$ of labelled $k$-gonal 2-trees on $n$ $k$-gons is
given by %
\begin{equation} %
\label{eq:aetiqkodd} 
a_n = {1\over 2}\left( m^{n-2}+1 \right),\quad \quad n\geq 2, 
\end{equation} %
where $m=(k-1)n+1$ is the number of edges.
\end{prop} 

\begin{figure}[h] 
  \centerline{\includegraphics[width=.8\textwidth]{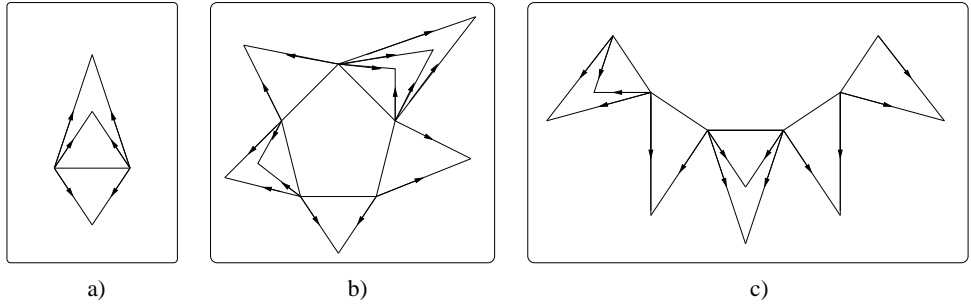}}
 \caption{Skeleton species a) $Q(X,Y)$, b) $S(X,Y)$ and c) $U(X,Y)$} 
 \label{fig:squelette} 
\end{figure}

For the unlabelled enumeration, 
notice from Figure~\ref{fig:page}a that in
every $k$-gon containing the pointed (but not oriented) edge of an
$\ab$-structure, it is possible to orient the $k-1$ other edges in a
canonical direction, "away from the root edge", when $k$ is odd (but there remains an
ambiguous opposite edge if $k$ is even). This phenomenon permits us to introduce
{\it skeleton} species, when $k$ is odd, in analogy with the approach
of Fowler et al. They are the two-sort
quotient species $Q(X,Y)$, $S(X,Y)$ and $U(X,Y)$, where $X$ represents
the sort of $k$-gons and $Y$ the sort of oriented edges, defined
by Figures~\ref{fig:squelette}a, b and c, where $k=5$.

In analogy with the case $k=3$, we get the following propositions.
\medskip
\begin{prop}
The skeleton species $Q$, $S$ and $U$ admit the following expressions
in terms of quotients species
\begin{eqnarray} 
Q(X,Y) &=& E(XY^2)/\Z_2,\\
S(X,Y) &=& C_k(E(XY^2))/\Z_2,\\ 
U(X,Y) &=& (E(XY^2))^k/\Z_2 .
\end{eqnarray} 
\end{prop}
\medskip
\begin{prop} %
\label{prop:pointed}
For $k$ odd, $k\geq 3$, we have the following expressions for the pointed 
species of $k$-gonal 2-trees, where $B=\arbre^{\rightarrow}$:
\begin{equation} 
\ab =   Q(X,B^{k-1\over 2}),\quad
\at =   X\cdot S(X,B^{k-1\over 2}) ,\quad 
\abt =  X\cdot U(X,B^{k-1\over 2}).
\end{equation} 
\end{prop}

In order to obtain enumerative formulas, we have to compute the
cycle index series of the species $Q$, $S$ and $U$.
\medskip
\begin{prop} %
\label{prop:ZQSU}
The cycle index series of the species $Q(X,Y)$, $S(X,Y)$ and
$U(X,Y)$ are given by
\begin{eqnarray}
Z_Q &=& {1\over 2}\bigg(Z_{E(XY^2)}+q\bigg),\label{eq:ZQ} \\ 
Z_S &=& {1\over 2}\bigg( Z_{C_k(E(XY^2))}+ %
  q\cdot (p_2\circ Z_{E(XY^2)})^{k-1\over 2}\bigg)\label{eq:ZS},\\ 
Z_U &=& {1\over 2}\bigg( Z_{(E(XY^2))^{k}}+q\cdot (p_2\circ Z_{E(XY^2)})^{k-1\over 2} \bigg),\label{eq:ZU} 
\end{eqnarray}  
where  %
\begin{equation} %
q=h\circ (x_1y_2+p_2\circ(x_1{y_1^2-y_2\over 2})),  %
\label{eq:q} %
\end{equation} %
$p_2$ represents the power sum symmetric function of degree two, $h$ the homogeneous
symmetric function and $\circ$, the plethystic substitution.
\end{prop}
\begin{preuve}  %
We use a two-sort extension of formula (\ref{eq:ZF/Z2}) but the sort $Y$ is the %
important one here. The variables ${\mathbf s}$ will keep track of the colored triangles %
and new variables ${\mathbf t}=(t_1,t_2,\dots)$, of the colored oriented edges  %
and we seek to express the inventory in terms %
of the power sums $x_i=p_i({\mathbf s})$ and $y_i=p_i({\mathbf t})$. %
Hence the second terms of the right-hand-sides %
of formulas (\ref{eq:ZQ})--(\ref{eq:ZU}), represent the $\tau$-symmetric unlabelled  %
colored $F(X,Y)$-structures.  %
For example, for (\ref{eq:ZQ}), the given formula (\ref{eq:q}) simply expresses the fact that  %
a $\tau$-symmetric unlabelled colored $Q(X,Y)$-structure consists of a set of \emph{pages},  %
where the $\tau$ symmetry comes either from a page with identically colored  %
oriented edges or from pairs of pages whose oriented edges are oppositely colored. %
See \cite{TF1,TF2} for more details. %

In the case of $S$, we have to leave fixed an unlabelled colored $C_k(E(XY^2))$-structure.  %
For this, the cycle of length $k$ must possess (at least) one symmetry axis
passing through the middle of one of its sides. The attached
structure on this distinguished edge must be globally left fixed;
this gives the factor $q$. On each side of the axis, each colored
$E(XY^2)$-structure must have its mirror image; this contributes 
the factor $(p_2\circ Z_{E(XY^2)})^{k-1\over 2}$. %
It can be seen that in the case of higher degree of symmetry, %
the choice of the symmetry axis is arbitrary.   %
The reasoning is very similar for the species $U$ and in fact the  %
$\tau$-symmetric term is the same as in the previous case.
\end{preuve} 

It is now a simple matter to combine the Dissymmetry Theorem with  %
Propositions \ref{prop:pointed} and \ref{prop:ZQSU} and the substitution rules of unlabelled
enumeration in order to obtain $\widetilde{\arbre}(x)$. %
Note that the first terms of formulas (\ref{eq:ZQ})--(\ref{eq:ZU}) %
will give rise to $\widetilde{\arbre}_o(x)$ and that a cancellation will %
occur in the $\tau$-symmetric terms, leaving only  %
$q(x_i\mapsto x^i,y_i\mapsto\widetilde{B}^{k-1\over2}(x^i))$ %
to compute.
\medskip
\begin{prop} 
Let $k\geq 3$ be an odd integer. The ordinary generating series
$\widetilde{\arbre}(x)$ of unlabelled $k$-gonal 2-trees is given by
\begin{equation} %
\widetilde{\arbre}(x) = \frac{1}{2}\bigg(\widetilde{\arbre}_o(x)+ %
\exp \big( \sum_{i\geq 1} \frac{1}{2i}(2x^i\widetilde{B}^{k-1\over2}(x^{2i}) + %
x^{2i}\widetilde{B}^{k-1}(x^{2i}) - x^{2i}\widetilde{B}^{k-1\over 2}(x^{4i})\big)\bigg). %
\label{ordinary} 
\end{equation}
\end{prop} 
\begin{cor} 
  For $k\geq 3$, odd, the number $\widetilde{a}_{n}$ of unlabelled
  $k$-gonal 2-trees over $n$ $k$-gons, satisfy the following
  recurrence
\begin{equation} 
\widetilde{a}_n=\frac{1}{2n}\sum_{j=1}^n\bigg(\sum_{l|j}l\omega_l\bigg)\Big(\widetilde{a}_{n-j}
-{1\over 2}\widetilde{a}_{o,n-j}\Big)+{1\over 2}\widetilde{a}_{o,n},\quad \widetilde{a}_0=1, 
\end{equation} 
where, for all $n\geq 1$,
\begin{equation} 
\omega_n=2b_{n-1\over 2}^{({k-1\over 2})}+b_{n-2\over 2}^{(k-1)}-b_{n-2\over 4}^{({k-1\over 2})}, 
\end{equation} 
and $b_i^{(j)}$ is defined in Corollary~\ref{corol}.
\end{cor} 
%
%
\subsection{Case $k$ even} 
The case $k$ even is more delicate. For example, as observed by one of the anonymous %
referees, there are more than one labelled oriented $k$-gonal 2-tree left fixed  %
by an orientation reversal. They can be obtained by taking an edge labelled %
ordinary tree and replacing edges by k-gons attached at opposite edges. %
These $k$-gonal 2-trees have no \emph{side decoration} and this explains %
their symmetry with respect to orientation.  %
It is known (and follows from (\ref{eq:betiqu}) for $k=2$)  %
that the number of edge-labelled trees with $n$ edges is $(n+1)^{n-2}$.  %
Hence we have the following: %
\begin{prop} \label{ankimpair}
If $k$ is even, 
the number $a_n$ of labelled $k$-gonal 2-trees on $n$ $k$-gons is
given by %
\begin{equation} %
\label{eq:aetiqkeven} 
a_n = {1\over 2}\left( m^{n-2} + (n+1)^{n-2} \right),\quad \quad n\geq 2, 
\end{equation} %
where $m=(k-1)n+1$ is the number of edges.
\end{prop}  %

For the unlabelled enumeration of the three species $\ab$,
$\at$ and $\abt$, we apply relation (\ref{eq:F/Z2tilde}) to 
formulas (\ref{eq:abarre}). For the species $\ab=\arbre^{\rightarrow}/\Z_2$,  %
the action of $\tau$ consists in reversing the orientation of the rooted edge. %
we have
\begin{equation} 
\widetilde{\arbre}^{-}(x)={1\over 2}(\widetilde{\arbre}^{\rightarrow}(x)
+\widetilde{\arbre}_{\tau}^{\rightarrow}(x)), 
\end{equation} 
where $\widetilde{\arbre}_{\tau}^{\rightarrow}(x)$ %
is the tilde generating series of $\tau$-symmetric (unlabelled)  %
oriented-edge-rooted 2-trees.  %
Let $\arbre_{\rm S}$ denote the subspecies of $B=\arbre^\rightarrow$
consisting of $\arbre^\rightarrow$-structures $s$ which are
isomorphic to their image $\tau\cdot s$. %
We have to compute $\widetilde{\arbre}_{\mathrm S}(x)=\widetilde{\arbre}_{\tau}^\rightarrow(x)$. %

\begin{figure}[h] 
\centerline{\includegraphics[width=.9\textwidth]{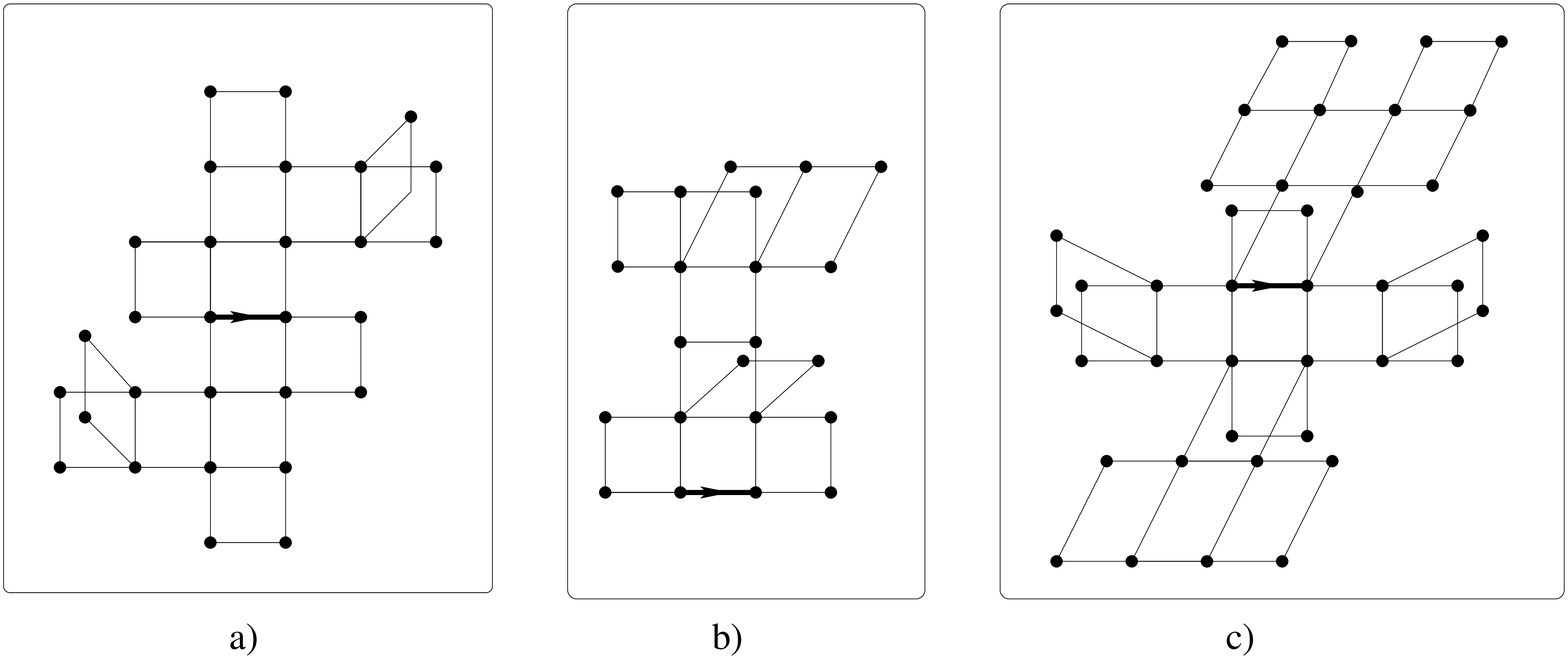}} 
\caption{a) A pair of alternated pages, b) a mixed page,  %
c) a totally symmetric $\arbre^\rightarrow$-structure} 
\label{fig:asnouveau} 
\end{figure}

Let us introduce some auxiliary subspecies of $\arbre_{\rm S}$ which appear %
when we analyse these $\tau$-symmetric structures in terms of their \emph{pages} %
that is their maximal sub-2-trees containing a unique triangle adjacent %
to the rooted edge. %
We say that there is some \emph{crossed symmetry}  %
if we can find, inside the 2-tree, two \emph{alternated} pages, that is pages of the form $\{s,\tau\cdot s\}$,  %
where $s$ is not itself $\tau$-symmetric, attached to the same root edge. %
See Figure \ref{fig:asnouveau}a %
Let $P_{\rm AL}$ denote the subspecies of pairs of alternated pages.  %
A {\it mixed page} is a symmetric page having at least one crossed symmetry.  %
See Figure~\ref{fig:asnouveau}b.  %
Let $P_{\rm M}$ denote the species of mixed pages. %

Finally, we say that a page is \emph{totally symmetric} or \emph{vertically symmetric}  %
if it contains no crossed symmetries.  %
Let $P_{\rm TS}$ denote the species of totally symmetric pages and  %
set  %
\begin{equation} %
\asym=E(P_{\rm TS}), \label{eq:aTS} %
\end{equation}  %
the subspecies of \emph{totally symmetric}  %
$\arbre^\rightarrow$-structures. See Figure~\ref{fig:asnouveau}c. %
We can characterize all these species and their tilde generating series %
by functional equations. First, we have
\begin{equation} %
P_{\rm TS} = X\cdot X_{=}^2<B^{k-2\over 2}>\cdot\asym, %
\label{eq:PTS}  
\end{equation}  %
where $X_{=}^2<F>$ represents the species of ordered pairs of isomorphic
$F$-structures.  %
Note that $(X_{=}^2<F>)^{\sim}(x)=\widetilde{F}(x^2)$. %
Translating equations (\ref{eq:aTS}) and (\ref{eq:PTS}) in terms of tilde generating
series, we get %
\begin{equation} %
\widetilde{\arbre}_{\rm TS}(x)  =   %
  \exp\left(\sum_{j\geq1}\widetilde{P}_{\rm TS}(x^j)\right) %
\label{eq:tildeaTS} %
\end{equation}  %
and %
\begin{equation} %
\widetilde{P}_{\rm TS}(x)  =  x\ \widetilde{B}^{k-2\over 2}(x^{2}) %
\widetilde{\arbre}_{\rm TS}(x). %
\label{eq:tildePTS} 
\end{equation} 
\begin{prop} \label{prop:pibeta}
The numbers $\pi_n=|\widetilde{P}_{\rm TS}[n]|$ and %
$\beta_n=|\widetilde{\arbre}_{\rm TS}[n]|$ of unlabelled totally symmetric pages  %
and $\arbre^\rightarrow$-structures, respectively, on $n$ polygons,  %
satisfy the following system of recurrences: $\beta_{0}=1$ and,  %
for $n\geq 1$, %
\begin{eqnarray} %
\pi_n &= & \sum_{i+j=n-1\atop i\ even}b_{i\over 2}^{({k-2\over 2})}\beta_j, %
\label{eq:PTSn} \\
\beta_n & = & {1\over n}\sum_{j=0}^{n-1}\beta_{j}\sum_{d|n-j}d\pi_d. %
\label{eq:asymn} %
\end{eqnarray}  %
\end{prop} 
\begin{preuve}  %
Formula (\ref{eq:PTSn}) is obvious. For (\ref{eq:asymn}), it suffices to take  %
$x$ times the logarithmic derivative of (\ref{eq:tildeaTS}).
\end{preuve}
%
%
%

Now, from the definition of the species $P_{\rm AL}$ of pairs of alternated pages, %
we have 
\begin{equation}
P_{\rm AL} = \Phi_{2}<XB^{k-1}-(P_{\rm TS}+P_{\rm M})>,\label{eq:PAL} %
\end{equation}
%
%
%
where $\Phi_{2}<F>$ represents the species of unordered pairs of $F$-structures
of the form $\{s,\tau\cdot s \}$. Note that  %
$\Phi_{2}<F>^{\sim}(x) = {1\over 2}\widetilde{F}(x^2)$ whenever the structures  %
$s$ and $\tau\cdot s$ are guaranteed not to be isomorphic, so that  %
\begin{equation}
\widetilde{P}_{\rm AL}(x) = {1\over 2} \big(x^2\widetilde{B}^{k-1}(x^2)  %
- \widetilde{P}_{\rm TS}(x^{2}) - \widetilde{P}_{\rm M}(x^{2})\big). %
\label{eq:PALtilde} %
\end{equation}
Also by definition, the 
species $P_{\rm M}$ of mixed pages satisfies %
\begin{eqnarray}
P_{\rm M} & = & X\cdot X_=^2<B^{k-2\over 2}>\cdot (\arbre_{\rm S}-\arbre_{\rm TS}) \nonumber\\ %
 & = & X\cdot X_=^2<B^{k-2\over 2}>\cdot \arbre_{\rm S} - P_{\rm TS}, %
\label{eq:pm} %
\end{eqnarray}
so that %
%
\begin{equation} 
\widetilde{P}_{\rm M}(x) 
= x\widetilde{B}^{k-2\over 2}(x^2)\widetilde{\arbre}_{\rm S}(x)-\widetilde{P}_{\rm TS}(x). %
\label{eq:PMtilde}
\end{equation}

\begin{figure}[h] 
\centerline{\includegraphics[width=.65\textwidth]{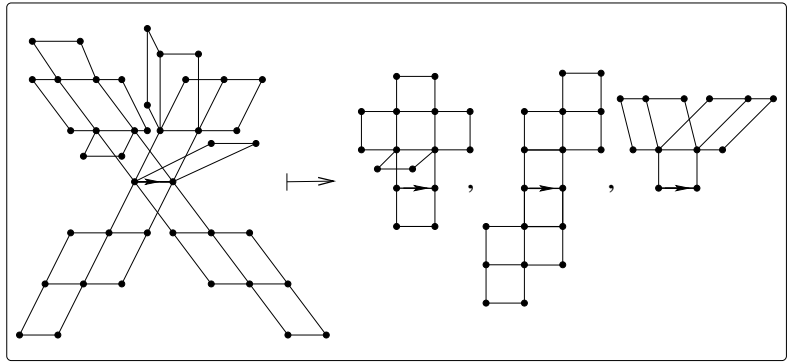}} 
\caption{Decomposition of a $\tau$-symmetric  $\widetilde{\arbre}^{\rightarrow}$-structure} \label{fixabarre} 
\end{figure}

Finally, for  the tilde generating series $\widetilde{\arbre}_{\rm S}(x)$  %
of unlabelled $\tau$-symmetric $\arbre^{\rightarrow}$-structures, we
have (see Figure~\ref{fixabarre})
\begin{eqnarray}
\widetilde{\arbre}_{\rm S} (x)& = & E(P_{\rm TS}+P_{\rm AL}+P_{\rm M})^{\sim}(x), %
\label{eq:aStilde}\\
 & = & \exp \bigg(\sum_{i\geq 1}{1\over i}( \widetilde{P}_{\rm TS}(x^i)+\widetilde{P}_{\rm AL}(x^i)
+\widetilde{P}_{\rm M}(x^i) )      \bigg).\label{eq:aStildebis}
\end{eqnarray} %
From equations (\ref{eq:PALtilde}), (\ref{eq:PMtilde}) and (\ref{eq:aStildebis}) %
we deduce the following. %
\begin{prop} %
The numbers $\alpha_n=\widetilde{a}_{{\rm S},n}$ of unlabelled $\tau$-symmetric  %
$\arbre^{\rightarrow}$-structures,
$\widetilde{P}_{{\rm AL},n}$, of pairs of alternated pages and  %
$\widetilde{P}_{{\rm M},n}$ of mixed pages, on $n$ $k$-gons are characterized by the %
following system of recurrences: $\alpha_{0}=1$, and for $n\geq1$, %
\begin{equation} 
\widetilde{P}_{{\rm M},n}=\sum_{i=0}^{n-1}b_{i\over 2}^{({k-2\over 2})}\alpha_{n-1-i}-\pi_{n},
\end{equation} %
\begin{equation} 
\widetilde{P}_{{\rm AL},n}={1\over 2}\left(b_{n-2\over 2}^{(k-1)}-\pi_{n/2}- %
\widetilde{P}_{{\rm M},n/2}\right),
\end{equation} %
\begin{equation} 
\alpha_{n} ={1\over n}\sum_{i=1}^{n}\bigg( \sum_{d|i}d\omega_d\bigg)\alpha_{n-i}, 
\end{equation} 
where $\pi_n=\widetilde{P}_{{\rm TS},n}$ is given by Propositon \ref{prop:pibeta} and %
\begin{equation}
\displaystyle{\omega_k=\pi_{k}+\widetilde{P}_{{\rm AL},k}+ \widetilde{P}_{{\rm M},k}}.
\end{equation} %
\end{prop}
\begin{prop} 
If $k$ is an even integer, then the number of unlabelled (unoriented) edge rooted
$k$-gonal 2-trees over $n$ $k$-gons is given by
\begin{equation} 
\widetilde{a}_n^- = {1\over 2}(b_n+\alpha_{n}).
\end{equation}
\end{prop} 
%
%
%
%

\begin{figure}[h] 
\centerline{\includegraphics[width=.40\textwidth]{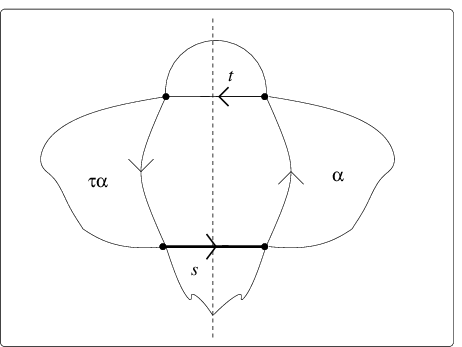}}
  \caption{A $\tau$-symmetric unlabelled $\arbre^{\btr}_o$-structures} 
  \label{fig:x1} 
\end{figure}

Let us now turn to the species $\abtpair$ of $k$-gonal 2-trees rooted
at an edge-pointed $k$-gon. 
\medskip
\begin{prop}
We have
\begin{equation} 
{\tildeabtpair}(x) = {1\over 2}\bigg(\tildeaobtpair(x)+ %
\widetilde{\arbre}_{o,\tau}^{\: \underline{\diamond\!}}(x)\bigg), 
\end{equation}
where   
\[
\widetilde{\arbre}_{o,\tau}^{\: \underline{\diamond\!}
}(x)=x\widetilde{\arbre}_{\rm S}^2(x)\widetilde{B}^{k-2\over 2}(x^2).
\]
\end{prop}
\begin{preuve}
An unlabelled $\tau$-symmetric $\aobt$-structure possesses an axis
of symmetry which is, in fact, the mediatrix of the distinguished
edge of the root polygon, and also the mediatrix of its opposite edge; %
see Figure~\ref{fig:x1}. The two structures
$s$ and $t$ glued on these two edges are thus symmetric, which leads
to the term $(\widetilde{\arbre}_{\rm S}(x))^2$. Then, on each side
of the axis, are found two $B^{k-2\over 2}$-structures $\alpha $ and
$\beta$, which by symmetry satisfy $\beta = \tau\cdot\alpha$,
contributing to the factor $\widetilde{B}^{k-2\over 2}(x^2)$.
\end{preuve}
\medskip
\begin{cor}
We have the following expression for the number
$\widetilde{a}_n^{{\: \underline{\!\diamond\!}}}$ of unlabelled
$\abtpair$-structures,
\begin{equation} 
\widetilde{a}_n^{{\: \underline{\diamond}} } = {1\over 2}\bigg(\widetilde{a}_{o,n}^{{\: \underline{\diamond}} } 
+\sum_{i+j=n-1}\alpha_{i}^{(2)}\cdot b_{j\over 2}^{({k-2 \over 2})}\bigg), 
\end{equation} 
where $\alpha_{i}^{(2)} =[x^i]\widetilde{\arbre}_{S}^2(x)$.\hfill$\Box$ 
\end{cor} 

We proceed in a similar way for the species $\atpair$, of $k$-gon
rooted $k$-gonal 2-trees. Once again, we use relation
(\ref{eq:F/Z2tilde}), giving
\begin{equation} 
\tildeatpair(x)={1\over 2}\bigg(\tildeaotpair(x)+ \widetilde{\arbre}_{o,\tau}^{\: {\diamond}}(x) \bigg). 
\end{equation} 
\medskip
\begin{prop}
Let $\widetilde{\arbre}_{o,\tau}^{\: \diamond}(x)$ be the generating
series of unlabelled $\arbre_{o}^{\: \diamond}$-structures which are left fixed
by orientation reversing. Then, we have
\begin{equation} \label{pastis}
\widetilde{\arbre}_{o,\tau}^{\: {\diamond}}(x)  = {x\over 2} \widetilde{\arbre}_{\rm S}^2(x)
\widetilde{B}^{k-2\over 2}(x^2)   + {x\over 2}\widetilde{B}^{k\over 2}(x^2).
\end{equation} 
\end{prop}
\begin{preuve}
Notice first that in order to be left fixed by orientation reversing, an
$\aot$-structure must admit a reflective symmetry, along an axis which
can either pass through the middle of two opposite edges, or
pass through opposite vertices of the pointed polygon.   %
The enumeration is carried out by first
orienting the axis of symmetry. The first term of (\ref{pastis}) then
corresponds to an edge--edge symmetry, and the second term to a vertex--vertex
symmetry. The structures having both symmetries are
precisely those which are counted one half time in both of these
terms.  This is established for a general $k$ by considering the
unique power of 2, $2^m$, such that $k/2^m$ is odd. We illustrate the
proof in the following lines with $k=12$; the reader will easily
convince himself of the validity of this argument for any $k$.
\begin{figure}[h] 
 \centerline{\includegraphics[width=.8\textwidth]{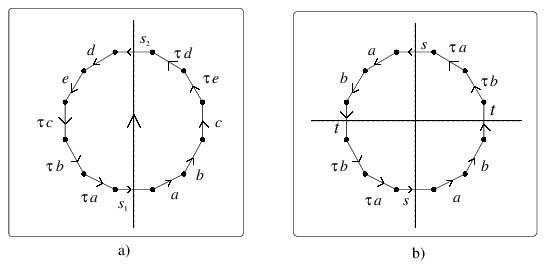}} 
  \caption{$\widetilde{\arbre}^{\tr}_{o,\tau}$-structures with an
    edge--edge symmetry}
  \label{fig:x2} 
\end{figure}

For $k=12$, a general unlabelled $\tau$-symmetric polygon-rooted
oriented $k$-gonal 2-tree with an oriented edge--edge axis will be of
the form illustrated in Figure~\ref{fig:x2} a), where $s_1$ and $s_2$
represent unlabelled $\arbre_{\rm S}$-structures, $a$, $b$, $c$, $d$
and $e$ are general unlabelled $B$-structures and $\tau x$
represents the opposite of the $B$-structures $x$, obtained by
reversing their orientation. Most of these structures are enumerated
exactly by ${1\over 2}x\widetilde{\arbre}_{\rm S}^2(x)\widetilde{B}^5(x^2)$. 
Indeed, the factor
$x\widetilde{\arbre}_{\rm S}^2(x)\widetilde{B}^5(x^2)$ is obtained in
the same way as for $\arbre_{o,\tau}^{\:\underline{\diamond}}$-structures and the division by two is
justified in the following cases:
\begin{enumerate}
\item $s_1\neq s_2$ (two orientations of the axis),
\item $s_1=s_2=s$, $(a,b,c)\neq (d,e,\tau\cdot c)$ (two orientations),
\item $s_1=s_2=s$, $(a,b,c)= (d,e,\tau\cdot c)$, so that $c=\tau\cdot c =t\in \widetilde{\arbre}_{\rm S}$, 
and either 
$s\neq t$\ or\ $s=t$ and $(a,b)\neq (\tau\cdot b,\tau\cdot
  a)$ (two choices for the symmetry axis, see Figure~\ref{fig:x2} b)),
\end{enumerate}
However, the structures with $s=t$ and $b=\tau\cdot a$ (see
Figure~\ref{fig:x3}) will occur only once and are counted only one
half time in the formula. But, notice that these structures also admit
a vertex--vertex symmetry axis and, as it will turn out, are also
counted one half time in the second term of (\ref{pastis}).
\begin{figure}[h] 
\centerline{\includegraphics[width=.35\textwidth]{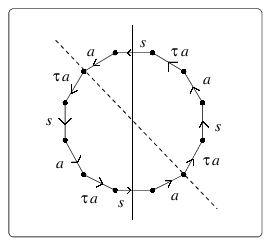}} 
\caption{An $\widetilde{\arbre}^{\tr}_{o,\tau}$-structure with
edge--edge and vertex--vertex symmetries}
\label{fig:x3} 
\end{figure}

Similarly, an unlabelled $\arbre_{o,\tau}^{\: \diamond}$-structure
with an oriented vertex--vertex symmetry axis will be of the form
illustrated in Figure~\ref{fig:x4} a), where $a$, $b$, $\ldots $, $f$
are arbitrary unlabelled $B$-structures. Most of these terms are
enumerated exactly by ${1\over 2}x\widetilde{B}^6(x^2)$, the division
by two being justified in the following cases:
\begin{enumerate}
\item $(a,b,c)\neq (d,e,f)$ (two orientations of the symmetry axis),
\item $(a,b,c)=(d,e,f)$ and $(a,b,c)\neq (\tau\cdot c,\tau\cdot
  b,\tau\cdot a)$ (two choices for the symmetry axis, see
  Figure~\ref{fig:x4} b)),
\end{enumerate}
\begin{figure}[h] 
 \centerline{\includegraphics[width=.8\textwidth]{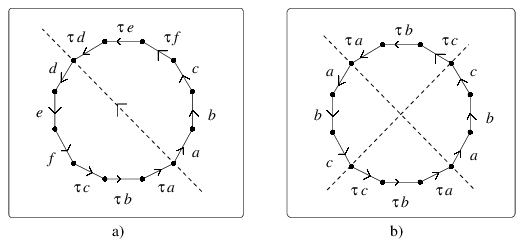}} 
  \caption{$\widetilde{\arbre}^{\tr}_{o,\tau}$-structures with a
  vertex--vertex symmetry axis}
  \label{fig:x4} 
\end{figure}
However, the structures with $(a,b,c)=(d,e,f)$, $c=\tau\cdot a$ and
$b=\tau\cdot b=s\in\widetilde{\arbre}_{\rm S}$ appear only once and
are counted one half time here. But they also have an edge-edge
symmetry axis and were also counted one half time in the first term of
(\ref{pastis}) (exchange $a$ and $\tau\cdot a$ in
Figure~\ref{fig:x3}).
\end{preuve}

The Dissymmetry Theorem yields, for $k$ even,
\begin{equation} 
\widetilde{\arbre}(x)  =  {1\over 2}\widetilde{\arbre}_o(x) %
 +{1\over 2}\widetilde{\arbre}_{\rm S}(x) %
  +{1\over 2}\widetilde{\arbre}_{o,\tau }^{\: \diamond }(x) %
  -{1\over 2}\widetilde{\arbre}_{o,\tau }^{\: \underline{\diamond} }(x),
\end{equation}
and we have the following result.
\medskip
\begin{prop}
Let $k$ be an even integer, $k\geq 4$. Then the generating series
$\widetilde{\arbre}(x)$ of unlabelled $k$-gonal 2-trees is given by
\begin{equation}
\widetilde{\arbre}(x) =  {1\over 2}\widetilde{\arbre}_o(x) %
+{1\over 2}\widetilde{\arbre}_{\rm S}(x) %
+{x\over 4}(\widetilde{B}^{k\over 2}(x^2) %
-\widetilde{\arbre}^2_{\rm S}(x)\widetilde{B}^{k-2\over 2}(x^2)). %
\label{ordinaryeven}
\end{equation}
\end{prop}
\begin{cor} 
Let $k$ be an even integer, $k\geq 4$. Then the number of unlabelled
$k$-gonal 2-trees over $n$ $k$-gons is given by
\begin{equation}\label{resultat1}
\widetilde{a}_{n} = {1\over 2}\widetilde{a}_{o,n} +{1\over 2}\alpha_{n} 
+ {1\over 4}b_{n-1\over 2}^{({k\over 2})} 
- {1\over 4}\sum_{i+j=n-1}\alpha_{i}^{(2)}\cdot b_{j\over 2}^{({k-2 \over 2})}, 
\end{equation} 
where
$$
b_{l}^{(m)}=[x^l]\widetilde{B}^m(x), \quad
\alpha_i^{(2)}=[x^i]\widetilde{\arbre}_{\rm S}^2(x).
$$  
\end{cor}
\medskip %

Note that the case $k=2$ corresponds to ordinary trees with $n$ edges  %
and that the formulas given here are also valid when properly interpreted. %
Table \ref{tableantilde} 
gives the exact values of the numbers
$\widetilde{a}_n$ of unlabelled $k$-gonal 2-trees with $n$ $k$-gons,  %
for $k$ from $2$ up to $12$ and for $n=0,1,\ldots,20$.
{
\begin{table}[h]
\footnotesize
{\boldmath$ k=2$}

1, 1, 1, 2, 3, 6, 11, 23, 47, 106, 235, 551, 1301, 3159, 7741, 19320, 
48629, 123867, 317955, 823065, 2144505


\noindent{\boldmath$ k=3$}

\noindent1, 1, 1, 2, 5, 12, 39, 136, 529, 2171, 9368, 41534, 188942, 874906, 
4115060, 19602156, 94419351, 459183768, 2252217207, 11130545494, 
55382155396


\noindent{\boldmath$ k=4$}

\noindent1, 1, 1, 3, 8, 32, 141, 749, 4304, 26492, 169263, 1115015, 7507211, 
51466500, 358100288, 2523472751, 17978488711, 129325796854, 
938234533024, 6858551493579, 50478955083341


\noindent{\boldmath$ k=5$}

\noindent1, 1, 1, 3, 11, 56, 359, 2597, 20386, 167819, 1429815, 12500748, 
111595289, 1013544057, 9340950309, 87176935700, 822559721606, 
7836316493485, 75293711520236, 728968295958626, 7105984356424859


\noindent{\boldmath$  k=6$}

1, 1, 1, 4, 16, 103, 799, 7286, 71094, 729974, 7743818, 84307887, 
937002302, 10595117272, 121568251909, 1412555701804, 16594126114458, 
196829590326284, 2354703777373055, 28385225424840078, 
344524656398655124


\noindent{\boldmath$ k=7$}

\noindent1, 1, 1, 4, 20, 158, 1539, 16970, 199879, 2460350, 31266165, 
407461893, 5420228329, 73352481577, 1007312969202, 14008437540003, 
196963172193733, 2796235114720116, 40038505601111596, 
577693117173844307, 8392528734991449808


\noindent{\boldmath$ k=8$}

\noindent1, 1, 1, 5, 26, 245, 2737, 35291, 483819, 6937913, 102666626, 
1558022255, 24133790815, 380320794122, 6081804068869, 98490990290897, 
1612634990857755, 26660840123167203, 444560998431678554, 
7469779489114328514, 126375763235359105446


\noindent{\boldmath$ k=9$}

1, 1, 1, 5, 32, 343, 4505, 66603, 1045335, 17115162, 289107854, 
5007144433, 88516438360, 1591949961503, 29053438148676, 
536972307386326, 10034276171127780, 189331187319203010, 
3603141751525175854, 69097496637591215442, 1334213677527481808220


\noindent{\boldmath$ k=10$}

1, 1, 1, 6, 39, 482, 7053, 117399, 2070289, 38097139, 723169329, 
14074851642, 279609377638, 5651139037570, 115901006038377, 
2407291353219949, 50553753543016719, 1071971262516091572, 
22926544048209731554, 494103705426160765546, 10722146465907412669810


\noindent{\boldmath$ k=11$}

\noindent1, 1, 1, 6, 46, 636, 10527, 194997, 3823327, 78118107, 1646300388, 
35570427615, 784467060622, 17601062294302, 400750115756742, 
9240636709048733, 215435023547580882, 5071520482516388865, 
120417032326341878672, 2881134828445365441407, 69410468220307148620226


\noindent{\boldmath$  k=12$}

\noindent1, 1, 1, 7, 55, 840, 15189, 309607, 6671842, 149850849, 3471296793, 
82442359291, 1998559329142, 49290785442796, 1233639304644946, 
31268489727956101, 801335133177932829, 20736286803363051714, 
541224489038545084067, 14234799536039481373552, 
376974819516101224941091
%
\caption{Values of $\widetilde{a}_{n}$ for $k=2,\ldots,12$ and $n=0,\ldots,20$}
\label{tableantilde}
\end{table}
}

\section{Asymptotics}
Thanks to the Dissymmetry Theorem and to the various combinatorial
equations related to it, the asymptotic enumeration of 
unlabelled $k$-gonal 2-trees depends essentially on the asymptotic
enumeration of $B$-structures where $B$ is the auxiliary species
characterized by the functional equation (\ref{eq:aoriente1}).  %

We first give the following result, which is a consequence of the
classical theorem of Bender (see \cite{Bender}) and is inspired from
the approach of Fowler et al. for 2-trees (see \cite{TF1,TF2}).
\medskip
\begin{prop}\label{growth} 
Let $p=k-1$. Let us write $b(x)=\widetilde{B}(x)=\sum b_n(p)x^n$.  %
Let $\xi_p$ be the smallest root of the equation
\begin{equation}\label{eq:xip}
\xi = {1\over ep}\omega^{-p}(\xi),
\end{equation}
where $\omega(x)$ is defined by
\begin{equation}
\omega(x)={\rm e}^{{1\over 2}x^2b^{p}(x^2)+{1\over 3}x^3b^p(x^3)+\cdots}.
\end{equation}
Then, there exist constants $\alpha_p$ and $\beta_p$ such that
\begin{equation}\label{eq:bnp}
b_n(p)\sim \alpha_p\beta_p^nn^{-{3/2}},\qquad \mbox{as } n\rightarrow \infty.
\end{equation}
Moreover, 
\begin{equation}\label{grille}
\alpha_p=\alpha(\xi_p)={1\over\sqrt{2\pi}}{1\over p^{1+{1\over p}}}\xi_p^{-{1\over p}}
\bigg(1+{p\xi_p\omega^{'}(\xi_p)\over\omega(\xi_p)}\bigg)^{1\over 2}
\end{equation} 
and
\begin{equation}\label{cercle}
\displaystyle{\beta_p={1\over \xi_p}},
\end{equation}
\end{prop}
\begin{preuve}
The functional equation (\ref{eq:btilde}) implies that $y=b(x)$ satisfies the relation %
\begin{equation} %
\label{eq:omega}
y=e^{xy^p}\omega(x).
\end{equation} %
By Bender's theorem applied to the function
$f(x,y)=y-e^{xy^p}\omega(x)$, we have to find a solution
$(\xi_p,\tau_p)$ of the system
\begin{equation}
f(x,y)=0 \qquad \mbox{and}\qquad f_y(x,y)=0.
\end{equation}
It is equivalent to say that $\xi_p$ is solution of (\ref{eq:xip}) and
that $p\xi_p\tau_p^p=1$. 
In fact, $\xi_p$ is the radius of convergence of $b(x) $ and $\sqrt{\xi_p}$ is
radius of convergence of $\omega(x)$. It can be
shown that $0<\xi_p<1$ so that $0<\xi_p<\sqrt{\xi_p}<1$. %
Indeed, if $\rho_p$ is the radius of convergence of the algebraic function $\theta(x)$  %
defined by $\theta=1+x\theta^p$, then, using Lagrange Inversion Formula and Stirling's Formula,  %
we obtain $\rho_p=(p-1)^{p-1}/p^p<1$, for $p\geq2$. Now, take a small fixed $x>0$  %
and consider the two curves  %
$z=\varphi_1(y)=1+xy^p$ and $z=\varphi_2=e^{xy^p}\omega(x)$ in the $(y,z)$-plane. %
Since $\varphi_1(y)<\varphi_2(y)$, for $y>0$, and $\theta(x)=\varphi_1(\theta(x))$ and %
$b(x)=\varphi_2(b(x))$, we have that $\theta(x)<b(x)$. If $x_0>\rho_p$, we must have  %
$b(x_0)=\infty$ since $\theta(x_0)=\infty$. This implies that $\xi_p\leq\rho_p$. %
For $p=1$ ($k=2$), a similar argument with $\varphi_1(y)=1+xy+xy^2/2$ shows that %
$\xi_1\leq\sqrt{2}-1$. %
Note also that from the recurrence (\ref{eq:bnonetiq}) it follows that $b_n(p)$ is %
bounded by the coefficient $c_n$ of the function $c(x)$ defined by $c=1+xc^k$, so that %
we have $\xi_p\geq\rho_{p+1}=p^p/(p+1)^{p+1}$, for $p\geq1$. %

Since $f_{yy}(\xi_p,\tau_p)\neq 0$, $\xi_p$ is an algebraic
singularity of degree 2 of $b(x)$ and, for $x$ near $\xi_p$, we have
an expression of the form
\begin{equation}\label{etoile}
b(x)=\tau_{p,0}+\tau_{p,1}(1-{x\over \xi_p})^{1\over 2}+\tau_{p,2}(1-{x\over \xi_p})
+\tau_{p,3}(1-{x\over \xi_p})^{3\over 2}+\cdots 
\end{equation}
where 
\begin{eqnarray}
\tau_{p,0} &=& \tau_p \ =\ b(\xi_p) \ = \ \left({1\over p\xi_p} \right)^{1\over p},\label{etoile1}\\
\tau_{p,1} &=& -{\sqrt{2}\over p^{1+{1\over p}}}\xi_p^{-{1\over p}}
\left( 1+\frac{p\xi_p\omega^{'}(\xi_p)}{\omega(\xi_p)} \right)^{1\over 2},\label{etoile2}\\
\tau_{p,2} &=& \frac{1}{3p^{2+{1\over p}}}\xi_p^{-{1\over p}}
\left( (2p+3)-p(p-3)\frac{\xi_p\omega^{'}(\xi_p)}{\omega(\xi_p)} \right).\label{etoile3}
\end{eqnarray}
The asymptotic formula (\ref{eq:bnp}) with $\alpha_p$ and $\beta_p$
given by (\ref{grille}) and (\ref{cercle}) then follow from the fact
that the main term of the asymptotic behavior of the coefficients
$b_n(p)$ of $x^n$ in (\ref{etoile}) depends only on the term
$\tau_{p,1}(1-{x\over \xi_p})^{1\over 2}$ in (\ref{etoile}) and is
given by
\begin{equation}\label{croixrond}
b_n(p) \sim  {{1\over 2 }\choose n}\tau_{p,1}(-1)^{n} {1\over \xi_p^{n}} 
\sim \alpha_p \beta_p^{n} n^{-{3\over 2}}\quad {\rm as} \quad n \rightarrow \infty .
\end{equation}
\end{preuve}
Note that numerical
approximations of $\xi_p$, for fixed $p$, can be computed by iteration
using $(\ref{eq:xip})$, and a suitable truncated polynomial
approximation of $b(x)$. We now state our main asymptotic result.
\medskip
\begin{prop}\label{bonhomme}
Let $p=k-1$. Then, the number $\widetilde{a}_n$ of $k$-gonal 2-trees
on $n$ unlabelled $k$-gons satisfy
\begin{equation}\label{diese}
\widetilde{a}_n \sim {1\over 2}\widetilde{a}_{o,n}, \quad n\rightarrow \infty,
\end{equation}
where $\widetilde{a}_{o,n}$ is the number of oriented $k$-gonal
2-trees over $n$ unlabelled polygons. Moreover, 
\begin{equation}\label{newgrid}
\widetilde{a}_{o,n} \sim \overline{\alpha}_p\beta_p^n n^{-{5/2}},
\quad n\rightarrow \infty,
\end{equation}
where 
\begin{eqnarray}
\overline{\alpha_p} &=& 2\pi p^{1+{2\over p}}\xi_p^{2\over p}\alpha_p^3,\label{alpha1}\\
 &=& {1\over \sqrt{2\pi}}{1\over p^{2+{1\over p}}}\xi_p^{-{1\over p}}
 \left(1+p\frac{\omega^{'}(\xi_p)}{\omega(\xi_p)}  \right)^{3\over 2}\label{carrediag},
\end{eqnarray}
and $\beta_p={1\over \xi_p}$ is the same growth as in
Proposition~\ref{growth}.
\end{prop}
\begin{preuve}
The asymptotic formula (\ref{diese}) follows from the fact that
the radius of convergence, $\xi_p$, of $\widetilde{\arbre}(x)$ 
is equal to the radius of convergence of the
dominating term ${1\over 2}\widetilde{\arbre}_{o}(x)$. This is due
to the easily checked fact that all terms in (\ref{ordinary}) and
(\ref{ordinaryeven}), except ${1\over 2}\widetilde{\arbre}_{o}(x)$,
have a radius of convergence greater or equal to
$\sqrt{\xi_p}>\xi_p$. To establish (\ref{newgrid}), note first that,
because of equation (\ref{eq:aoriente}), the radius of convergence
of $\widetilde{\arbre}_o(x)$ is equal to the radius of convergence,
$\xi_p$, of
\begin{equation}\label{spiraldream}
b(x)-{k-1\over k}xb^k(x),
\end{equation} 
where $b(x)=\widetilde{B}(x)$ and $k=p+1$. This implies that the
asymptotic behavior of the coefficients $\widetilde{a}_{o,n}$ of
$\widetilde{\arbre}_o(x)$ is completely determined by that of
(\ref{spiraldream}). Substituting (\ref{etoile}) into (\ref{spiraldream}) 
and making use of (\ref{etoile3}) gives the following expansion 
\begin{equation}\label{carrepoint}
b(x)-{k-1\over k}xb^k(x) = \overline{\tau}_{p,0}+\overline{\tau}_{p,1}
\left(1-{x\over \xi_p}  \right)^{1\over 2}+\overline{\tau}_{p,2}\left(1-{x\over \xi_p}  \right)
+\overline{\tau}_{p,3}\left(1-{x\over \xi_p}  \right)^{3\over 2}+\cdots 
\end{equation}
where 
\begin{eqnarray}
\overline{\tau}_{p,0} &=& {p\over p+1}\tau_{p,0},\label{0barre}\\
\overline{\tau}_{p,1} &=& 0,\label{1barre}\\
\overline{\tau}_{p,2} &=& -{1\over 2}\frac{p(p+1)\tau^2_{p,1}-2\tau^2_{p,0}}{(p+1)\tau_{p,0}},\label{2barre}\\
\overline{\tau}_{p,3} &=& -{1\over 6}\frac{\tau_{p,1}(6p\tau_{p,0}\tau_{p,2}+
p(p-1)\tau^2_{p,1}-6\tau^2_{p,0} )}{\tau^2_{p,0}},\label{3barre}\\
   &=&  -{p\over 3}{\tau^3_{p,1}\over \tau^2_{p,0}}.\label{4barre}
\end{eqnarray}
This implies that the dominating term for the asymptotic behavior of 
the coefficients $\widetilde{a}_{n,o}$ of $x^n$ in $\widetilde{\arbre}_o(x)$ 
depends only on the term 
$\overline{\tau}_{p,3}\left(1-{x\over \xi_p}  \right)^{3\over 2}$ in (\ref{carrepoint}) 
and is given by 
\begin{equation}\label{grospoint}
\widetilde{a}_{n,o} \sim {{3\over 2}\choose n}\overline{\tau}_{p,3}(-1)^n{1\over \xi_p^n}
\sim\overline{\alpha}_p\beta_pn^{-{5\over 2}},\quad \mbox{as } n\rightarrow \infty.
\end{equation}
Computations making use of (\ref{4barre}), (\ref{etoile1}) and
(\ref{etoile2}), show that $\overline{\alpha}_p$ is indeed given by
(\ref{alpha1}) and (\ref{carrediag}).
\end{preuve}

Our final result gives an explicit formula in terms of integer
partitions for the common radius of convergence $\xi_p$ of the series
$\widetilde{B}(x)$, $\widetilde{\arbre}(x)$ and
$\widetilde{\arbre}_o(x)$ from which the growth constant
$\beta_p={1\over \xi_p}$ is obtained. We need the following special
notations. If $\lambda=(\lambda_1\geq\lambda_2\geq \ldots \geq
\lambda_{\nu})$ is a partition of an integer $n$ in $\nu$ parts, we
write $\lambda \vdash n$, $n=|\lambda|$, $\nu=l(\lambda)$,
$m_i(\lambda)=|\{j:\lambda_j=i\}|$ = number of parts of size $i$ in
$\lambda$. Furthermore, we put
\begin{equation}
\sigma_i(\lambda)=\sum_{d|i}dm_d(\lambda), \quad \sigma^*_i(\lambda)=\sum_{d|i,d<i}dm_d(\lambda),
\end{equation}
\begin{equation}
\widehat{\lambda}=1+|\lambda|+l(\lambda),\quad \widehat{z}(\lambda)
=2^{m_1(\lambda)}m_1(\lambda)!3^{m_2(\lambda)}m_2(\lambda)!\ldots.
\end{equation}
\medskip
\begin{prop}
  We have the convergent expansion
\begin{equation}\label{eq:xipsomme}
 \xi_p=\sum_{n=1}^{\infty}{c_n\over p^n},
\end{equation}
 where the coefficients $c_n$ are constants, independent of $p$,
 explicitely given by
\begin{equation}\label{eq:cn}
 c_n=\sum_{\lambda\vdash n}{e^{-\widehat{\lambda}}\over \widehat{\lambda}\widehat{z}(\lambda)}
 \prod_{i\geq 1}(\sigma_i(\lambda)-\widehat{\lambda})^{m_i(\lambda)-1} (\sigma_i^*(\lambda)-\widehat{\lambda}), 
\end{equation}
 where $\lambda$ runs over the set of partitions of $n$.
\end{prop}
\begin{preuve}
  We establish the explicit formulas (\ref{eq:xipsomme}) and
  (\ref{eq:cn}) by applying first Lagrange inversion to the equation
  $\xi=zR(\xi)$ where $z={1\over {\rm e}p}$ and $R(t)=\omega^{-p}(t)$,
  to get
\begin{equation}\label{pasbeau}
\xi_p = \xi =\sum_{n\geq 1}\gamma_n\left({1\over {\rm e}p } \right)^n, 
\quad {\rm and}\quad \gamma_n={1\over n}[t^{n-1}]\omega^{-np}(t).
\end{equation}
Next, to explicitely evaluate $\omega^{-np}(x)$, we use Labelle's version
(\cite{1234}) of the Good inversion formula in the context of cycle
index series as follows. We begin with
\begin{eqnarray}
\omega^p(x) &=& \exp({1\over 2}px^2b^{p}(x^2)+{1\over 3}px^3b^p(x^3)+\cdots),\\
            &=& \exp({1\over 2}px_2+{1\over 3}px_3+\cdots )\circ Z_{XB^p(X)}\bigg|_{x_i:=x^i}
            \label{pointdexclamation}
\end{eqnarray}
where the $\circ $ denotes the plethystic substitution. Using
(\ref{cest5}), we can then write $XB^p(X)={A(pX)\over p}$. This
implies that
\begin{equation}\label{rondcroix}
\omega^p(x) = \exp({1\over 2}px_2+{1\over 3}px_3+\cdots )\circ \frac{Z_{A}(px_1,px_2,\ldots)}{p}\bigg|_{x_i:=x^i},  
\end{equation}
and we get
\begin{eqnarray}
\omega^{-np}(x) &=& \exp(-{n\over 2}px_2-{n\over 3}px_3-\cdots)
\circ\left( {1\over p}Z_A(px_1,px_2,\ldots ) \right)\bigg|_{x_i:=x^i} \\
            &=& \exp(-{n\over 2}x_2-{n\over 3}x_3-\cdots)\circ Z_A(x_1,x_2,\ldots )\bigg|_{x_i:=px^i}.
            \label{carrepetitecroix} 
\end{eqnarray}
Then, using Labelle's inversion formula for cycle index series, we have, for
any formal cycle index series $g(x_1,x_2,\ldots )$
\begin{equation}\label{saletete}
[x_1^{n_1}x_2^{n_2}\ldots ]\ g\circ Z_A(x_1,x_2,\ldots ) = [t_1^{n_1}t_2^{n_2}\ldots ]g(t_1,t_2,\ldots )
\prod_{i=1}^{\infty}(1-t_i)\exp(n_i(t_i+{1\over 2}t_{2i}+\cdots    )),
\end{equation}
and 
\begin{equation}\label{bonnetete}
\prod_{j=1}^{\infty}\exp(n_j(t_j+{1\over 2}t_{2j}+\cdots )) = \prod_{i=1}^{\infty}\exp(\sum_{d|i}dn_d{t_i\over i}).
\end{equation}
Taking $g(x_1,x_2,\ldots )= \exp(-{\nu\over 2}px_2-{\nu\over 3}px_3-\cdots)$, 
gives, after some computations,
\[
[x_1^{n_1}x_2^{n_2}\ldots ]\left( \exp(-{\nu\over 2}x_2-{\nu\over 3}x_3-\cdots)\circ Z_{A}  \right) = \hspace{6cm}
\]
\begin{equation}\label{oeil}
\left\{ 
\begin{array}{ccc}
0 & {\rm if} & n_1>0,\\ 
\left(\frac{\displaystyle{\prod_{i\geq2}}(-\nu+\displaystyle{\sum_{d|i}}dn_d)^{n_i-1}
(-\nu+\displaystyle{\sum_{d|i,d<i}}dn_d)}{\displaystyle{2^{n_2}n_2!3^{n_3}n_3!\ldots}}\right)
& {\rm if} & n_1=0.\\
\end{array}
\right.
\end{equation}
Making the substitution $x_i:=px^{i}$, for $i=1,2,3,\ldots $, gives the 
explicit formula
\[
\displaystyle{\omega^{-\nu p}(x) = \sum_{n\geq 0}\left(\sum_{2n_2+3n_3+\cdots =n}p^{n_2+n_3+\cdots } 
\frac{{\displaystyle{\prod_{i\geq 2}}}(-\nu+\displaystyle{\sum_{d|i}}dn_d)^{n_i-1}
(-\nu+\displaystyle{\sum_{d|i,d<i}}dn_d)}{2^{n_2}n_2!3^{n_3}n_3!\ldots }      \right)x^{n}}.
\]
This implies, taking $\nu=n$ and using (\ref{pasbeau}), that
\begin{eqnarray*}
\xi_p &=& \displaystyle{\sum_{n\geq 1}{1\over n}\left(\sum_{2n_2+3n_3+\cdots =n-1}p^{n_2+n_3+\cdots }
\frac{\displaystyle{ \prod_{i\geq 2}}(1-n+\displaystyle{\sum_{d|i}}dn_d)^{n_i-1}
(1-n+\displaystyle{\sum_{d|i,d<i}}dn_d)}{2^{n_2}n_2!3^{n_3}n_3!\ldots }      \right)
\left( {1\over {\rm e}p}\right)^{n}}, \\
      &=& \sum_{n\geq 1}{c_n\over p^{n}},
\end{eqnarray*}
where the coefficients $c_n$, $n\geq 1$, are given by (\ref{eq:cn}).
\end{preuve}

Here are the first few values of the universal constants $c_n$
occuring in (\ref{eq:xipsomme}), for $n=1,\ldots ,5$.
\begin{eqnarray}
c_1 & = & {1\over {\rm e}}\ =\  0.36787944117144232160,\nonumber \\
c_2 & = & -{1\over 2}{1\over {\rm e}^3}\ =\  -0.02489353418393197149,\nonumber \\
c_3 & = & {1\over 8}{1\over {\rm e}^5}-{1\over 3}{1\over {\rm e}^4} \ =\  -0.00526296958802571004, \\
c_4 & = & -{1\over 48}{1\over {\rm e}^7}+{1\over {\rm e}^6}-{1\over 4}{1\over {\rm e}^5}\ =\ 0.00077526788594593923,
\nonumber\\
c_5 & = & {1\over 384}{1\over {\rm e}^9}-{4\over 3}{1\over {\rm e}^8}+{49\over 72}{1\over {\rm e}^7}
-{1\over 5}{1\over {\rm e}^6}\ =\ 0.00032212622183609932.\nonumber
\end{eqnarray}
%

Table \ref{tableconstants} 
gives, to 12 decimal places, the
constants $\xi_p$, $\alpha_p$, $\overline{\alpha}_p$ and
$\beta_p={1\over \xi_p}$ for $p=1,\ldots ,12$.  %
%

\begin{table}[h]
\small
\centerline{ 
\begin{tabular}{|c|c|c|c|c|}
\hline 
$p$ & $\xi_p$ & $\alpha_p$ & $\overline{\alpha}_p$   &$\beta_p$ \\
\hline 
1 & \  0.338321856899\ 
  &\   1.300312124682\ 
  &\   1.581185475409\
  &\   2.955765285652\ 
\\ 
2 &\   0.177099522303 \ 
  &\   0.349261381742\ 
  &\   0.349261381742
  &\   5.646542616233\  
\\
3 &\  0.119674100436\   
  &\  0.191997258650\
  &\  0.067390781222\
  &\  8.356026879296\
\\
4 &\  0.090334539604\   
  &\  0.131073637349\   
  &\  0.034020667269\
  &\ 11.069962877759\
\\
5 &\  0.072539192528\    
  &\  0.099178841365\    
  &\  0.020427915489\
  &\ 13.785651110085\
\\
6 &\  0.060597948397\   
  &\  0.079660456931\   
  &\  0.013601784466\
  &\ 16.502208844693\
\\
7 &\  0.052031135998\       
  &\  0.066517090385\    
  &\  0.009699566188\
  &\ 19.219261329064\
\\
8 &\  0.045585869619\    
  &\  0.057075912245\    
  &\  0.007262873797\
  &\ 21.936622211299\
\\
9 &\  0.040561059517\    
  &\  0.049970993036\    
  &\  0.005640546218\
  &\ 24.654188324989\
\\
10&\  0.036533820306\   
  &\  0.044433135893\   
  &\  0.004506504206\
  &\ 27.371897918664\
\\
11&\  0.033233950789\     
  &\  0.039996691773\  
  &\  0.003682863427\
  &\ 30.089711763681\
\\ 
\hline 
\end{tabular} 
}
\caption{Numerical values of $\xi_p$, $\alpha_p$, $\overline{\alpha}_p$
and $\beta_p$, $p=1,\ldots,12$}
\label{tableconstants}
\end{table}
%
%

\begin{rem}
The computations of this section are also valid for the case $k=2$ 
($p=1$), corresponding to the class of ordinary rooted trees ({\it Cayley
trees}) defined by the functional equation $A=XE(A)$. In this case, %
the growth constant $\beta=\beta_1$, in {\rm (\ref{eq:bnp})}, is
known as the Otter constant (see {\rm \cite{Otter}}). It is
interesting to note that this constant takes the explicit form  %
$\beta={1\over \xi_1}$, with %
\begin{equation}
\xi _1= \sum_{n\geq 1}c_n.
\end{equation}
\end{rem}

\section*{Acknowledgments} %
We thank the referees for correcting a mistake in Proposition \ref{ankimpair} %
and for making many constructive suggestions.
%
%

 
%
{\it E-mail adresses}: [gilbert,leroux]@lacim.uqam.ca, lamathe@loria.fr


%
\end{document}